# Neural Network Optimal Feedback Control with Guaranteed Local Stability

**Tenavi Nakamura-Zimmerer[1][2], Qi Gong[1], Wei Kang[3][1]**

[1]Department of Applied Mathematics, Baskin School of Engineering, University of California, Santa Cruz, CA
[2]Flight Dynamics Branch, NASA Langley Research Center, Hampton, VA. This work was performed while affiliated with the University of California.
[3]Department of Applied Mathematics, Naval Postgraduate School, Monterey, CA
CORRESPONDING AUTHOR: T. Nakamura-Zimmerer (e-mail: tenakamu@ucsc.edu)

This work was supported with funding from the Air Force Office of Scientific Research (AFOSR) under grant FA9550-21-1-0113, the National Science Foundation (NSF) under grant nos. 2134235 and 2202668, and the University of California, Santa Cruz, Baskin School of Engineering Dissertation Year Fellowship.

**ABSTRACT** Recent research shows that supervised learning can be an effective tool for designing near-optimal feedback controllers for high-dimensional nonlinear dynamic systems. But the behavior of neural network controllers is still not well understood. In particular, some neural networks with high test accuracy can fail to even locally stabilize the dynamic system. To address this challenge we propose several novel neural network architectures, which we show guarantee local asymptotic stability while retaining the approximation capacity to learn the optimal feedback policy semi-globally. The proposed architectures are compared against standard neural network feedback controllers through numerical simulations of two high-dimensional nonlinear optimal control problems: stabilization of an unstable Burgers-type partial differential equation, and altitude and course tracking for an unmanned aerial vehicle. The simulations demonstrate that standard neural networks can fail to stabilize the dynamics even when trained well, while the proposed architectures are always at least locally stabilizing. Moreover, the proposed controllers are found to be close to optimal in testing.

## I. INTRODUCTION

Designing optimal feedback controllers for high-dimensional nonlinear systems remains an outstanding challenge for the control community. Even when the system dynamics are known, to design such controllers one needs to solve a Hamilton-Jacobi-Bellman (HJB) partial differential equation (PDE), whose dimension is the same as that of the state space. This leads to the well-known curse of dimensionality, which rules out traditional discretization-based approaches.

Recent work has demonstrated the promise of supervised learning as one potential approach for handling such challenging, high-dimensional problems. The main idea is to generate data by solving many open loop optimal control problems (OCPs) and then fit a model to this data set, thus obtaining an approximate optimal feedback controller. Various specific model design and training approaches have been developed within this framework. Earlier work [1]–[3] uses sparse grid interpolation to approximate the solution of the HJB equation – called the *value function* – and its gradient, which is used to compute the optimal feedback control. This line of work has been futher developed using neural networks (NNs) [4]–[8] and sparse polynomials [9], significantly increasing the maximum feasible problem dimension. Alternatively, one can directly approximate the value gradient [10], [11] or control policy [4], [5], [11]–[14].

There are also several closely-related research directions which can be classified as self-supervised learning methods. The method of successive approximations is a well-studied approach based on iterative updates of a value function model and/or control policy by approximately solving a series of Lyapunov equations [15]–[17]. These methods are equipped with convergence guarantees but they often depend on specific problem dynamics, a priori access to a semi-globally stabilizing controller, or polynomial model structures whose size grow exponentially with the problem dimension. A second group of methods attempt to solve the HJB PDE in the least-squares sense by minimizing the residual of the PDE and boundary conditions at (randomly sampled) collocation points [18]–[21]. More recently, [22]–[24] have proposed methods to solve the HJB equation along its characteristics without generating data. Generally speaking, such self-supervised approaches avoid the cost of data generation by taking on a harder learning problem.

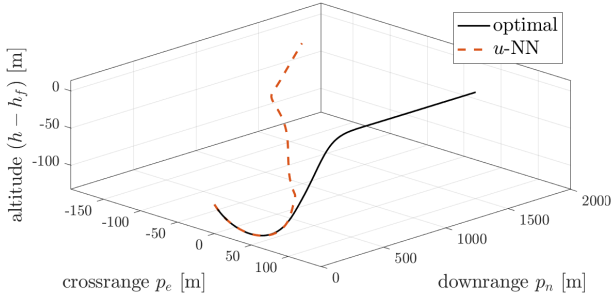

Figure 1: Simulated trajectory of a fixed wing UAV (29) controlled with an NN feedback controller, compared to the optimal trajectory.

Despite promising developments in the methodology, much less work has been done to study and improve the stability and reliability of these NN controllers. To see why this is needed, if we train a set of NNs to control a fixed wing unmanned aerial vehicle (UAV) studied in Section VI, a surprisingly large fraction of these fail to stabilize the system despite having good approximation accuracy. Figure 1 shows a closed loop simulation with one such controller where the NN-controlled trajectory closely tracks the optimal (stable) trajectory before suddenly destabilizing and eventually settling at an undesired steady state. Behavior like this obviates the need for better understanding, more rigorous testing, and more reliable algorithms.

This problem has been previously recognized by [11], [25], while [13] point out that test accuracy incompletely characterizes NN controller performance, suggesting some practical evaluations of optimality and stability. [26] study linear stability near a desired equilibrium, linear time delay stability, and stability around a nominal trajectory using high order Taylor maps. In terms of algorithm development, [8], [11] propose several NN architectures incorporating a linear quadratic regulator (LQR) which makes training more reliable and improves local stability properties.

The purpose of this paper is twofold: first, to bring attention to stability issues with NN-controlled systems; and second, to propose some NN architectures which can mitigate some of these challenges. These architectures improve on previous work in [8], [11] by guaranteeing, at a minimum, local asymptotic stability (LAS) of the system. This is accomplished by exactly recovering the LQR gain at the origin by construction. We also prove a universal approximation theorem for NNs with such structures, showing that they can approximate the nonlinear optimal feedback law up to arbitrary accuracy, and consequently provide semi-global stability and optimality.

This paper is organized as follows. We start by describing the problem setting in Section II. In Section III we describe the novel NN architectures and present theory underpinning their local stability properties and their ability to approximate the full nonlinear optimal feedback policy. In Section IV we outline the supervised learning procedure. It should be noted, however, that the proposed architectures do not have to be trained using supervised learning; in principle they can also be implemented with self-supervised learning. In Section V we apply several practical closed loop stability and optimality tests to demonstrate the advantages of the proposed NN architectures. As a testbed we use the Burgers'-type PDE system (25), which is nonlinear, open loop unstable, and high-dimensional. Then in Section VI we apply the proposed control design methodology to learn optimal feedback controllers for a six degree-of-freedom (6DoF) fixed wing UAV with nonlinear dynamics and aerodynamics, showing that the framework can be applied to practical problems. A summary and directions for future work are given in Section VII. The code used for simulations in this paper will be made publicly available at https://github.com/Tenavi/QRnet.

## II. PROBLEM SETTING

We focus on infinite horizon nonlinear OCPs of the form

$$\begin{cases} \underset{\mathbf{u}(\cdot)}{\text{minimize}} & J\left[\mathbf{u}(\cdot); \mathbf{x}_0\right] = \int_0^\infty \mathcal{L}(\mathbf{x}, \mathbf{u}) dt, \\ \text{subject to} & \dot{\mathbf{x}}(t) = \mathbf{f}(\mathbf{x}, \mathbf{u}), \\ & \mathbf{x}(0) = \mathbf{x}_0, \\ & \mathbf{u}(t) \in \mathbb{U}. \end{cases} \quad (1)$$

Here $\mathbf{x} : [0, \infty) \to \mathbb{R}^n$ is the state, $\mathbf{u} : [0, \infty) \to \mathbb{U} \subseteq \mathbb{R}^m$ is the control, $\mathbf{f} : \mathbb{R}^n \times \mathbb{U} \to \mathbb{R}^n$ is a vector field which is continuously differentiable ($\mathcal{C}^1$) in $\mathbf{x}$ and $\mathbf{u}$, and $(\mathbf{x}_f, \mathbf{u}_f) \in \mathbb{R}^n \times \mathbb{U}$ is a (possibly unstable) equilibrium of $\mathbf{f}(\cdot)$. We consider box control constraints

$$\mathbb{U} = \{\mathbf{u} \in \mathbb{R}^m \,|\, u_{\min,i} \leq u_i \leq u_{\max,i}, i = 1, \ldots, m\}, \quad (2)$$

for vectors $\mathbf{u}_{\max}, \mathbf{u}_{\min} \in \mathbb{R}^m$; and running costs $\mathcal{L} : \mathbb{R}^n \times \mathbb{U} \to [0, \infty)$ of the form

$$\mathcal{L}(\mathbf{x}, \mathbf{u}) = q(\mathbf{x}) + r(\mathbf{u}), \quad (3)$$

for smooth functions $q : \mathbb{R}^n \to [0, \infty)$, $r : \mathbb{U} \to [0, \infty)$ satisfying $q(\mathbf{x}_f) = 0$, $r(\mathbf{u}_f) = 0$, $q(\mathbf{x}) \geq 0$ for $\mathbf{x} \neq \mathbf{x}_f$, and $r(\mathbf{u}) > 0$ for $\mathbf{u} \neq \mathbf{u}_f$. This is a standard running cost for regularization or set-point tracking problems. We make the standard assumptions that $\mathbf{u}_f$ is an interior point of $\mathbb{U}$ and that the OCP (1) is well-posed, i.e. there exists an optimal control $\mathbf{u}^* : [0, \infty) \to \mathbb{U}$ such that $J[\mathbf{u}^*(\cdot)] < \infty$.

Due to real-time application requirements, we would like to design a control policy in explicit feedback form, $\mathbf{u} = \mathbf{u}^*(\mathbf{x})$, which can be evaluated online given any measurement of $\mathbf{x}$. The mathematical framework for designing such an optimal feedback policy is the HJB equation.

Define the *value function* $V : \mathbb{R}^n \to \mathbb{R}$ as the optimal cost-to-go of (1) starting at $\mathbf{x}(0) = \mathbf{x}$, i.e. $V(\mathbf{x}) \coloneqq J[\mathbf{u}^*(\cdot); \mathbf{x}]$. Under appropriate conditions, the value function is the unique viscosity solution [27] of the steady state HJB PDE,

$$\min_{\mathbf{u} \in \mathbb{U}} \mathcal{H}(\mathbf{x}, V_\mathbf{x}, \mathbf{u}) = 0, \qquad V(\mathbf{x}_f) = 0, \quad (4)$$

where $V_\mathbf{x} \coloneqq [\partial V / \partial \mathbf{x}]^T$ and

$$\mathcal{H}(\mathbf{x}, \boldsymbol{\lambda}, \mathbf{u}) \coloneqq \mathcal{L}(\mathbf{x}, \mathbf{u}) + \boldsymbol{\lambda}^T \mathbf{f}(\mathbf{x}, \mathbf{u}) \quad (5)$$

is the Hamiltonian. If (4) can be solved (in the viscosity sense), then it provides both necessary and sufficient conditions for optimality. Furthermore, the optimal feedback control is then obtained from the Hamiltonian minimization condition,

$$\mathbf{u}^*(\mathbf{x}) = \mathbf{u}^*(\mathbf{x}; V_\mathbf{x}(\mathbf{x})) = \arg\min_{\mathbf{u} \in \mathbb{U}} \mathcal{H}(\mathbf{x}, V_\mathbf{x}, \mathbf{u}). \quad (6)$$

## III. ARCHITECTURES FOR OPTIMAL FEEDBACK DESIGN

Our goal is to construct a feedback policy which approximates the optimal control, i.e. $\widehat{\mathbf{u}}(\mathbf{x}) \approx \mathbf{u}^*(\mathbf{x})$. While previous work has clearly demonstrated the potential of deep learning as a means of overcoming the curse of dimensionality in optimal control, NNs are notoriously "black boxes" and their behavior – especially when implemented in the closed loop system – is complex and hard to predict. Notably, even if we can train a highly accurate NN, it can still fail to stabilize the system. Thus there is a clear need for designing NN feedback controllers with *built-in* stability properties.

In prior work, the authors proposed $V$-*QRnet* [8] (originally just called *QRnet*), $\lambda$-*QRnet*, and $u$-*QRnet* [11]. These architectures combine an LQR controller with NNs. The LQR terms are good approximations of the optimal control near $\mathbf{x}_f$, and improve local stability. Meanwhile, the NNs are intended to capture nonlinearities and thereby learn the nonlinear optimal feedback over a large domain.

However, none of these architectures *guarantee* LAS, which motivates us to pursue alternative designs. In this paper we introduce four novel NN architectures, $\lambda_{\text{Jac}}$-*QRnet*, $\lambda_{\text{mat}}$-*QRnet*, $u_{\text{Jac}}$-*QRnet*, and $u_{\text{mat}}$-*QRnet*, all of which guarantee LAS of $\mathbf{x}_f$ while retaining the ability to approximate the nonlinear optimal control semi-globally.

The remainder of this section is organized as follows. In Section A we review $\lambda$-*QRnet*, $u$-*QRnet*, and LQR control design. The novel NN architectures are presented in Sections B and C. In Section D we show that these controllers automatically provide LAS, and finally in Section E we prove that they have the capacity to approximate the nonlinear optimal control.

Throughout the paper, architectures denoted with a leading $V$ approximate the value function, those with $\lambda$ approximate the value gradient, and those with $u$ directly approximate the optimal control. $V$-NN, $\lambda$-NN, and $u$-NN refer to standard feedforward NNs for approximating the value function, value gradient, and optimal control, respectively.

### A. $\lambda$-QRnet AND u-QRnet

Introduced in preliminary work [11], $\lambda$-*QRnet* and $u$-*QRnet* are straightforward linear combinations of LQR with NNs. Let

$$\begin{cases} \mathbf{A} \coloneqq \frac{\partial \mathbf{f}}{\partial \mathbf{x}}(\mathbf{x}_f, \mathbf{u}_f), & \mathbf{B} \coloneqq \frac{\partial \mathbf{f}}{\partial \mathbf{u}}(\mathbf{x}_f, \mathbf{u}_f), \\ \mathbf{Q} \coloneqq \frac{\partial^2 q}{\partial \mathbf{x}^2}(\mathbf{x}_f), & \mathbf{R} \coloneqq \frac{\partial^2 r}{\partial \mathbf{u}^2}(\mathbf{u}_f). \end{cases} \quad (7)$$

Under the standard conditions that $(\mathbf{A}, \mathbf{B})$ is controllable and $(\mathbf{A}, \mathbf{Q}^{1/2})$ is observable, LQR gives a quadratic value function approximation,

$$V^{\text{LQR}}(\mathbf{x}) = (\mathbf{x} - \mathbf{x}_f)^T \mathbf{P}(\mathbf{x} - \mathbf{x}_f), \quad (8)$$

and a linear state feedback law,

$$\mathbf{u}^{\text{LQR}}(\mathbf{x}) = \mathbf{u}_f - \mathbf{K}(\mathbf{x} - \mathbf{x}_f), \quad \mathbf{K} = \mathbf{R}^{-1}\mathbf{B}^T\mathbf{P}, \quad (9)$$

where $\mathbf{P} \in \mathbb{R}^{n \times n}$ is a positive definite matrix satisfying the Riccati equation,

$$\mathbf{Q} + \mathbf{A}^T\mathbf{P} + \mathbf{P}\mathbf{A} - \mathbf{P}\mathbf{B}\mathbf{R}^{-1}\mathbf{B}^T\mathbf{P} = \mathbf{0}. \quad (10)$$

Sufficiently near $\mathbf{x}_f$, the LQR value function (8) and control (9) are good approximations of the true value function $V(\cdot)$ and optimal control $\mathbf{u}^*(\cdot)$, respectively. Specifically, $[\partial \mathbf{u}^*/\partial \mathbf{x}](\mathbf{x}_f) = -\mathbf{K}$ and $[\partial^2 V/\partial \mathbf{x}^2](\mathbf{x}_f) = 2\mathbf{P}$. But further away from $\mathbf{x}_f$, the control is suboptimal and in some cases may even fail to stabilize the nonlinear dynamics. For this reason we are interested in combining the local stability and optimality of LQR with NNs to learn the full *nonlinear* optimal control $\mathbf{u}^*(\cdot)$ over a semi-global domain.

Now we describe $\lambda$-*QRnet*, which approximates the value gradient $V_\mathbf{x}(\cdot)$ as

$$\widehat{\boldsymbol{\lambda}}(\mathbf{x}) = 2\mathbf{P}(\mathbf{x} - \mathbf{x}_f) + \mathcal{N}(\mathbf{x}; \boldsymbol{\theta}) - \mathcal{N}(\mathbf{x}_f; \boldsymbol{\theta}). \quad (11)$$

Here $\mathcal{N} : \mathbb{R}^n \times \mathbb{R}^p \to \mathbb{R}^n$ is an NN with $\mathcal{C}^1$ activation functions and parameters $\boldsymbol{\theta} \in \mathbb{R}^p$, and the linear component $2\mathbf{P}(\mathbf{x} - \mathbf{x}_f)$ is the gradient of the LQR value function (8). We then substitute (11) into (6) to obtain an approximate optimal feedback control:

$$\widehat{\mathbf{u}}(\mathbf{x}) = \mathbf{u}^*\left(\mathbf{x}; \widehat{\boldsymbol{\lambda}}(\mathbf{x})\right). \quad (12)$$

$\lambda$-*QRnet* can be easily implemented when we can solve (6) for an explicit formula the optimal feedback control in terms of the state and value gradient, as is the case for many problems of interest.

Alternatively, we can directly approximate the optimal control with $u$-*QRnet*:

$$\widehat{\mathbf{u}}(\mathbf{x}) = \sigma\left[\text{sat}\left(\mathbf{u}^{\text{LQR}}(\mathbf{x})\right) + \mathcal{N}(\mathbf{x}; \boldsymbol{\theta}) - \mathcal{N}(\mathbf{x}_f; \boldsymbol{\theta})\right], \quad (13)$$

where now $\mathcal{N} : \mathbb{R}^n \times \mathbb{R}^p \to \mathbb{R}^m$. In (13), $\text{sat}(\cdot)$ is the saturation function defined for each $i = 1, \ldots, m$ as

$$[\text{sat}(\mathbf{u})]_i \coloneqq \begin{cases} u_{\min,i}, & u_i < u_{\min,i}, \\ u_i, & u_{\min,i} \le u_i \le u_{\max,i}, \\ u_{\max,i}, & u_{\max,i} < u_i. \end{cases} \quad (14)$$

Next, $\sigma : \mathbb{R}^m \to \mathbb{U}$ is a generalized logistic function which smoothly saturates the nonlinear control[1]:

$$\sigma(\mathbf{u}) := \mathbf{u}_{\min} + \frac{\mathbf{u}_{\max} - \mathbf{u}_{\min}}{1 + \mathbf{c}_1 \exp\left[-\mathbf{c}_2 \left(\mathbf{u} - \mathbf{u}_f\right)\right]}. \quad (15)$$

Here multiplication and division are performed element-wise, and we set the constants $\mathbf{c}_1, \mathbf{c}_2 \in \mathbb{R}^m$ as

$$\mathbf{c}_1 = \frac{\mathbf{u}_{\max} - \mathbf{u}_f}{\mathbf{u}_f - \mathbf{u}_{\min}}, \quad \mathbf{c}_2 = \frac{\mathbf{u}_{\max} - \mathbf{u}_{\min}}{(\mathbf{u}_{\max} - \mathbf{u}_f)(\mathbf{u}_f - \mathbf{u}_{\min})}. \quad (16)$$

It is straightforward to verify that these choices of $\mathbf{c}_1, \mathbf{c}_2$ satisfy $\sigma(\mathbf{u}_f) = \mathbf{u}_f$ and $[\partial \sigma / \partial \mathbf{u}](\mathbf{u}_f) = 1$. Consequently, $\sigma(\cdot)$ smoothly imposes saturation constraints while preserving the unsaturated behavior near $\mathbf{u}_f$, as we visualize in Figure 2.

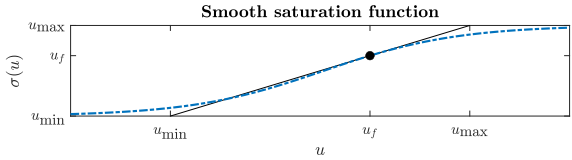

Figure 2: The smooth saturation function (15).

Use of this smooth saturation function makes learning easier since it prevents vanishing gradients when the control becomes saturated. The hard saturation function (14) acting on the LQR control inside the smooth saturation function (15) may appear redundant, but this is actually important for limiting the effect of the linear term away from $\mathbf{x}_f$. This reduces the burden on the NN: it is free to approximate nonlinearities without negating any excess contributions from the linear component.

It is easy to show that subtracting $\mathcal{N}(\mathbf{x}_f; \boldsymbol{\theta})$ in (11) and (13) makes the goal state $\mathbf{x}_f$ a closed loop equilibrium [11]. This is *not* true for standard NN controllers, $V$-NN, $\lambda$-NN, $u$-NN, or $V$-*QRnet*. Including the LQR terms improve local stability due to LQR's large gain and phase margins, but does not exactly recover LQR and so cannot assure LAS without adequate training.

### B. "JACOBIAN" QRnet ARCHITECTURES

Now we describe $\lambda_{\text{Jac}}$-*QRnet* and $u_{\text{Jac}}$-*QRnet*. These are similar to $\lambda$-*QRnet* and $u$-*QRnet*, except that we subtract the Jacobian of the NN components. This ensures that the controllers exactly recover LQR at $\mathbf{x}_f$, thus guaranteeing LAS. Furthermore, we will prove that these architectures retain the nonlinear function approximation capacity of standard feedforward NNs, allowing them to approximate the full nonlinear value gradient and optimal control.

First we have $\lambda_{\text{Jac}}$-*QRnet*:

$$\widehat{\boldsymbol{\lambda}}(\mathbf{x}) = \left[2\mathbf{P} - \frac{\partial \mathcal{N}}{\partial \mathbf{x}}(\mathbf{x}_f; \boldsymbol{\theta})\right](\mathbf{x} - \mathbf{x}_f) + \mathcal{N}(\mathbf{x}; \boldsymbol{\theta}) - \mathcal{N}(\mathbf{x}_f; \boldsymbol{\theta}). \quad (17)$$

---
[1]We also clip the exponent $-\mathbf{c}_2(\mathbf{u} - \mathbf{u}_f)$ to prevent numerical overflow when evaluating the gradient during training.

$u_{\text{Jac}}$-*QRnet* has an analogous structure:

$$\widehat{\mathbf{u}}(\mathbf{x}) = \sigma\left[\operatorname{sat}\left(\mathbf{u}^{\text{LQR}}(\mathbf{x})\right) - \left[\frac{\partial \mathcal{N}}{\partial \mathbf{x}}(\mathbf{x}_f; \boldsymbol{\theta})\right](\mathbf{x} - \mathbf{x}_f) + \mathcal{N}(\mathbf{x}; \boldsymbol{\theta}) - \mathcal{N}(\mathbf{x}_f; \boldsymbol{\theta})\right]. \quad (18)$$

These models are slower to train than $\lambda$-*QRnet* (11) and $u$-*QRnet* (13) since the Jacobian $[\partial \mathcal{N}/\partial \mathbf{x}](\mathbf{x}_f; \boldsymbol{\theta})$ must be evaluated during each forward pass. *After training, however, we can store the Jacobian matrix in memory so that it does not have to be recomputed online.* Therefore online evaluation is just as fast as $\lambda$-*QRnet* and $u$-*QRnet*.

### C. "MATRIX" QRnet ARCHITECTURES

In this section we describe $\lambda_{\text{mat}}$-*QRnet* and $u_{\text{mat}}$-*QRnet*. These alternatives to the "Jacobian"-style architectures employ *matrix-valued* NNs. Thus they avoid the costly Jacobian computations in exchange for having to optimize more NN parameters. These "matrix" *QRnet*s enjoy the same stability and approximation properties as the "Jacobian" *QRnet*s. We have not found a consistent performance advantage of either the "Jacobian" or "matrix" *QRnet*s: their relative learning ability appears to be problem-dependent.

First consider $\lambda_{\text{mat}}$-*QRnet*:

$$\widehat{\boldsymbol{\lambda}}(\mathbf{x}) = [2\mathbf{P} + \mathcal{N}(\mathbf{x}; \boldsymbol{\theta}) - \mathcal{N}(\mathbf{x}_f; \boldsymbol{\theta})](\mathbf{x} - \mathbf{x}_f). \quad (19)$$

Notice that in this case $\mathcal{N} : \mathbb{R}^n \times \mathbb{R}^p \to \mathbb{R}^{n \times n}$ is matrix-valued. Next we have $u_{\text{mat}}$-*QRnet*:

$$\widehat{\mathbf{u}}(\mathbf{x}) = \sigma\big[\operatorname{sat}\left(\mathbf{u}^{\text{LQR}}(\mathbf{x})\right) + [\mathcal{N}(\mathbf{x}; \boldsymbol{\theta}) - \mathcal{N}(\mathbf{x}_f; \boldsymbol{\theta})](\mathbf{x} - \mathbf{x}_f)\big], \quad (20)$$

where now $\mathcal{N} : \mathbb{R}^n \times \mathbb{R}^p \to \mathbb{R}^{m \times n}$.

A drawback of $\lambda_{\text{mat}}$-*QRnet* (19) is that the number of NN parameters scales with $O\left(Lw^2 + wn^2\right)$, where $L$ is the number of layers and $w$ is their width. For high dimensional OCPs, this can make (19) challenging to train as well as deploy on small processors. Meanwhile, the number of NN parameters in (20) scales with $O\left(Lw^2 + wmn\right)$. Since we typically have $m \ll n$, $u_{\text{mat}}$-*QRnet* is often much smaller and hence much faster to train than $\lambda_{\text{mat}}$-*QRnet*.

### D. LOCAL ASYMPTOTIC STABILITY GUARANTEES

Like $\lambda$-*QRnet* and $u$-*QRnet*, the new architectures automatically make the goal state $\mathbf{x}_f$ an equilibrium. Moreover, if we linearize the feedback control $\widehat{\mathbf{u}}(\cdot)$ at $\mathbf{x}_f$ then we recover the LQR control gain (9). This holds even when the models are poorly trained. This is desirable because LQR locally asymptotically stabilizes $\mathbf{x}_f$, and hence the proposed controllers are locally stabilizing by construction. This property is stated formally in Proposition 1 below. The proof is straightforward but tedious, so we omit it for brevity.

**Proposition 1 (Local asymptotic stability):**
Suppose $\widehat{\mathbf{u}}(\cdot)$ is a feedback policy specified by (12) with (17) or (19); or by (18) or (20). Then $[\partial \widehat{\mathbf{u}}/\partial \mathbf{x}](\mathbf{x}_f) = -\mathbf{K}$ and $\mathbf{x}_f$ is a locally stable equilibrium of the NN-controlled system, $\dot{\mathbf{x}} = \mathbf{f}(\mathbf{x}, \widehat{\mathbf{u}}(\mathbf{x}))$.

### E. APPROXIMATION CAPACITY

LAS is a critical but bare minimum requirement. To achieve the ultimate goal of semi-global stability and optimality through training, the NN architectures must be able to approximate $\mathbf{u}^*(\cdot)$ with sufficient accuracy. Unfortunately, we cannot directly use the Taylor series-like forms of (17–20) or existing NN universal approximation theorems like [28] to show this is possible. This is because $V_\mathbf{x}(\cdot)$ and $\mathbf{u}^*(\cdot)$ are in general *not* $\mathcal{C}^1$ everywhere, and because the NN architectures used in this work are not standard.

Nevertheless, for OCPs like (1) we expect $V_\mathbf{x}(\cdot)$ and $\mathbf{u}^*(\cdot)$ to be everywhere continuous and *locally* $\mathcal{C}^1$ in a neighborhood of $\mathbf{x}_f$. In this case Theorems 1 and 2 presented below show that NNs of the form (17–20) are universal approximators for such functions[2] To prove Theorems 1 and 2 we will first specialize the Stone-Weierstrass approximation theorem [29] to locally $\mathcal{C}^1$ functions, and then apply an NN universal approximation theorem [28]. For clarity of presentation, we simply state the main results here and defer the proofs to the Appendix.

Throughout this section let $\mathbb{X} \subset \mathbb{R}^n$ be compact, let $\mathbf{x}_f$ be an interior point of $\mathbb{X}$, and without loss of generality let $\mathbf{x}_f = \mathbf{0}$ and $\mathbf{u}_f = \mathbf{0}$. By $\mathcal{C}(\mathbb{X}; \mathbb{R}^d)$ we denote the space of continuous functions on $\mathbb{X}$ taking values in $\mathbb{R}^d$.

Our first main result concerns the approximation capacity of the "Jacobian" *QRnet* architectures introduced in Section B. As mentioned previously, since we expect the value gradient and optimal control for the OCP (1) to be continuous and locally $\mathcal{C}^1$, this supports the use of $\lambda_\text{Jac}$-*QRnet* and $u_\text{Jac}$-*QRnet* as nonlinear function approximators.

**Theorem 1 (Jacobian *QRnet* approximation):**
Suppose $\mathbf{f} \in \mathcal{C}(\mathbb{X}; \mathbb{R}^d)$, $\mathbf{f}(\mathbf{0}) = \mathbf{0}$, and $\mathbf{f}(\cdot)$ is $\mathcal{C}^1$ in a neighborhood of $\mathbf{0}$. Then for all $\varepsilon > 0$, there exists a feedforward NN with bounded, non-constant, $\mathcal{C}^1$ activation functions, $\mathcal{N} \in \mathcal{C}^1(\mathbb{X}; \mathbb{R}^d)$, such that for all $\mathbf{x} \in \mathbb{X}$,

$$\left\| \mathbf{f}(\mathbf{x}) - \left( \left[ \tfrac{\partial \mathbf{f}}{\partial \mathbf{x}}(\mathbf{0}) - \tfrac{\partial \mathcal{N}}{\partial \mathbf{x}}(\mathbf{0}) \right] \mathbf{x} + \mathcal{N}(\mathbf{x}) - \mathcal{N}(\mathbf{0}) \right) \right\|_1 < \varepsilon.$$

An analogous approximation theorem can be obtained for the "matrix" *QRnet* architectures introduced in Section C, $\lambda_\text{mat}$-*QRnet* and $u_\text{mat}$-*QRnet*.

**Theorem 2 (Matrix *QRnet* approximation):**
Suppose $\mathbf{f} \in \mathcal{C}(\mathbb{X}; \mathbb{R}^d)$, $\mathbf{f}(\mathbf{0}) = \mathbf{0}$, and $\mathbf{f}(\cdot)$ is $\mathcal{C}^1$ in a neighborhood of $\mathbf{0}$. Then for all $\varepsilon > 0$, there exists a feedforward NN with bounded, non-constant, $\mathcal{C}^1$ activation functions, $\mathcal{N} \in \mathcal{C}^1(\mathbb{X}; \mathbb{R}^{d \times n})$, such that for all $\mathbf{x} \in \mathbb{X}$,

$$\left\| \mathbf{f}(\mathbf{x}) - \left[ \tfrac{\partial \mathbf{f}}{\partial \mathbf{x}}(\mathbf{0}) + \mathcal{N}(\mathbf{x}) - \mathcal{N}(\mathbf{0}) \right] \mathbf{x} \right\| < \varepsilon.$$

---

[2] For simplicity, Theorems 1 and 2 do not address the saturation constraints (15) which may be applied to $u_\text{Jac}$-*QRnet* and $u_\text{mat}$-*QRnet*. In practice we find that the smooth saturation function does not hinder learning.

## IV. MODEL TRAINING

In this section we provide an overview of the supervised learning approach we use to train the NN controllers. Note that the proposed NNs do *not* have to be trained using supervised learning: they can be implemented in conjunction with *any* learning algorithm as long as an LQR controller can be computed for the system. In this work we focus on the impact of NN architecture rather than learning algorithm, so we restrict the scope to supervised learning.

Supervised learning can be broken down into three steps: data generation (Section A), NN optimization (Section B), and finally model evaluation against test data (Section C). In Section V we will illustrate a more rigorous test regimen specifically for control design, by which we compare the proposed controllers with LQR and standard NN controllers.

### A. DATA GENERATION

To circumvent the challenge of directly solving the HJB PDE (4), we can find an approximate optimal control by joining the solutions of many open loop OCPs. Each open loop OCP can be solved independently without the use of a spatial grid, thus mitigating the curse of dimensionality. Numerical methods based on this idea are referred to as *causality-free* [1].

To generate training and testing data sets for supervised learning we solve the open loop OCP (1) for a set of (randomly sampled) initial conditions. Note that in practice we approximate (1) by a finite horizon problem with large final time. Each open loop optimal trajectory and control provide input-output pairs $\mathbf{x}^{(i)}, \left( V_\mathbf{x}\left(\mathbf{x}^{(i)}\right), \mathbf{u}^*\left(\mathbf{x}^{(i)}\right) \right)$, where the superscript $(i)$ is the sample index. Aggregating data from all open loop solutions[3], we obtain a data set

$$\mathcal{D}_\text{train} = \left\{ \mathbf{x}^{(i)}, V_\mathbf{x}\left(\mathbf{x}^{(i)}\right), \mathbf{u}^*\left(\mathbf{x}^{(i)}\right) \right\}_{i=1}^{N_\text{train}}. \quad (21)$$

In the following we briefly review common computational methods for solving open loop optimal control. For more detailed discussions of data generation approaches for supervised learning, we refer the reader to [30] and references therein.

Algorithms for solving the open loop OCP (1) can be broadly classified as *indirect* and *direct* methods [31]. Indirect methods take the "optimize then discretize" approach, computing open loop optimal solutions to (1) by numerically solving necessary conditions of optimality from Pontryagin's Minimum Principle (PMP). These necessary conditions are given in term of a two-point boundary value problem (BVP) in terms of the state $\mathbf{x} : [0, \infty) \to \mathbb{R}^n$ and *costate* $\boldsymbol{\lambda} : [0, \infty) \to \mathbb{R}^n$, which under some conditions is equivalent to the value gradient along the optimal trajectory. There are mature BVP solvers that can be used for such computations.

Direct methods take the "discretize then optimize" approach, transforming the OCP (1) into a nonlinear programming problem. In contrast to indirect methods, direct methods do not require deriving the costate dynamics or initial guesses for the costates, and can more easily handle complicated OCPs such as those with path constraints. In the context of supervised learning, [13], [14] use a Hermite-Simpson method to generate data for finite horizon OCPs. Radau pseudospectral collocation [32], [33] is a direct method which is ideal for solving infinite horizon open loop OCPs. Pseudospectral methods have the added benefit of the covector mapping theorem [33], [34], which allows one to extract costate data from the solution of the discretized OCP.

---

[3] Note that there is no need to distinguish data from different trajectories as the value function and optimal feedback control are time-independent.

In this paper we employ an indirect method for the Burgers'-type PDE stabilization problem (Section V) and a Radau direct method for the UAV problem (Section VI).

### B. SUPERVISED LEARNING

Once a set of training data is available, the next step is training – i.e. data-driven optimization. Denoting the model parameters (i.e. the NN weights and biases) by $\boldsymbol{\theta} \in \mathbb{R}^p$, then the NN is trained by minimizing a mean squared error (MSE) loss function:

$$\boldsymbol{\theta} = \arg \min_{\boldsymbol{\theta}} \frac{1}{N_{\text{train}}} \sum_{i=1}^{N_{\text{train}}} \left\| \widehat{\mathbf{u}}\left(\mathbf{x}^{(i)}; \boldsymbol{\theta}\right) - \mathbf{u}^*\left(\mathbf{x}^{(i)}\right) \right\|_2^2. \quad (22)$$

As is standard in machine learning, the models learn on data which has been scaled to the range $[-1, 1]$, and the output is accordingly rescaled to the original domain when ultimately used for control.

When training the value gradient models, one can augment the loss function (22) with an additional MSE term to learn the value gradient,

$$\text{loss}_{\boldsymbol{\lambda}}(\boldsymbol{\theta}) = \frac{1}{N_{\text{train}}} \sum_{i=1}^{N_{\text{train}}} \left\| \widehat{\boldsymbol{\lambda}}\left(\mathbf{x}^{(i)}; \boldsymbol{\theta}\right) - V_{\mathbf{x}}\left(\mathbf{x}^{(i)}\right) \right\|_2^2, \quad (23)$$

and/or a term to minimize the residual of the HJB equation (4). The proposed NNs would also work well in conjunction with active learning methods [7].

### C. QUANTIFYING MODEL ACCURACY

To quantify the accuracy of the model we generate a second test data set, $\mathcal{D}_{\text{test}}$, from *independently drawn* initial conditions. During training, the NN sees only data points from the training set $\mathcal{D}_{\text{train}}$, while $\mathcal{D}_{\text{test}}$ is reserved for evaluating approximation accuracy after training. A typical metric is the relative mean $\ell^2$ error,

$$\text{RM}\ell^2 := \frac{\frac{1}{N_{\text{test}}} \sum_{i=1}^{N_{\text{test}}} \left\| \widehat{\mathbf{u}}\left(\mathbf{x}^{(i)}\right) - \mathbf{u}^*\left(\mathbf{x}^{(i)}\right) \right\|_2}{\max_{i=1,\ldots N_{\text{test}}} \left\| \mathbf{u}^*\left(\mathbf{x}^{(i)}\right) \right\|_2}, \quad (24)$$

where $N_{\text{test}}$ denotes the number of test points $\mathbf{x}^{(i)} \in \mathcal{D}_{\text{test}}$. A low test error indicates that the NN generalizes well, i.e. it did not overfit the training data.

However, even with low test error, there is a chance that the NN could still perform poorly when implemented in the closed loop system as seen in Figure 1. For this reason we believe that test metrics like (24) are insufficient in the context of control design; we should instead focus on rigorous closed loop stability and performance tests such as those presented in Section V.

## V. NUMERICAL RESULTS

In this section we compare the proposed controllers to standard feedforward NNs trained to approximate the value function, its gradient, and the optimal control. We also compare to $V$-QRnet [8], $\lambda$-QRnet, and $u$-QRnet [11]. We present results for three different tests:

1) linear stability near $\mathbf{x}_f$ (Section B);
2) Monte Carlo (MC) nonlinear stability (Section C);
3) MC optimality analysis (Section D).

Such tests are of course familiar to the control community, but we believe it is worth emphasizing their importance since more rigorous and realistic testing is needed in order to start trusting NN controllers in real-world applications. We also note that these tests are just a starting point: further examples include stabilization time [13], time delay stability [26], and robustness to measurement noise, disturbances, and parameter variations.

The numerical results clearly illustrate that standard NNs are *not* consistently stable, even when they have good approximation accuracy. Meanwhile, the results confirm that the proposed architectures guarantee LAS while still being able to accurately approximate the nonlinear optimal control throughout the training domain.

### A. UNSTABLE BURGERS'-TYPE PDE

To test the NN architectures we revisit the Burgers'-type PDE stabilization OCP from [8], [11]. This is a high-dimensional nonlinear OCP formulated by Chebyshev pseudospectral spatial discretization of an unstable version of a Burgers' PDE. Briefly, the problem can be summarized as

$$\begin{cases} \min_{\mathbf{u}(\cdot)} & J\left[\mathbf{u}(\cdot)\right] = \int_0^\infty \left(\mathbf{x}^T \mathbf{Q} \mathbf{x} + \mathbf{u}^T \mathbf{R} \mathbf{u}\right) dt, \\ \text{s.t.} & \dot{\mathbf{x}} = -\frac{1}{2}\mathbf{D}\mathbf{x} \circ \mathbf{x} + \nu \mathbf{D}^2 \mathbf{x} + \boldsymbol{\alpha} \circ \mathbf{x} \circ e^{-\beta \mathbf{x}} + \mathbf{B}\mathbf{u}. \end{cases} \quad (25)$$

Here $\mathbf{x} : [0, \infty) \to \mathbb{R}^n$ represents the PDE state $X(t, \xi)$ collocated at spatial coordinates $\xi_j = \cos(j\pi/n)$, $j = 1, \ldots, n$, $\mathbf{u} : [0, \infty) \to \mathbb{R}^m$ is the control, $\mathbf{D} \in \mathbb{R}^{n \times n}$ is the Chebyshev differentiation matrix, $\mathbf{Q} \in \mathbb{R}^{n \times n}$, $\mathbf{R} \in \mathbb{R}^{m \times m}$ are diagonal positive definite matrices, and "$\circ$" denotes element-wise multiplication. The parameters $\nu, \beta > 0, \boldsymbol{\alpha} \in \mathbb{R}^n$, and $\mathbf{B} \in \mathbb{R}^{n \times m}$ are defined in [8], and we take $n = 64$ and $m = 2$. Initial conditions are also selected as in [8].

We generate data by solving the OCP (1) for randomly sampled initial conditions, using an indirect method. For this problem we reliably obtain high quality data with the *SciPy* [35] implementation of the two-point BVP solver [36]. To get models with varying approximation accuracy, we generate training data sets with different numbers of trajectories.

Note that because data generation depends on random sampling and (22) is a highly non-convex optimization problem, results can vary for different random seeds. To account for this, for each different data set size we conduct ten trials with different randomly generated training trajectories and NN weight initializations. We evaluate the RM$\ell^2$ error (24) on an independent test data set containing 500 trajectories.

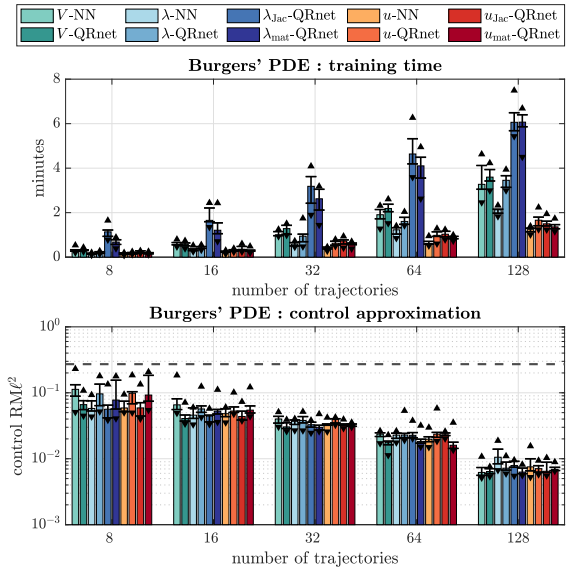

Figure 3: Training time and RM$\ell^2$ error for different NN architectures, depending on the amount of training data. Bar heights show the medians over ten trials, error bars show the 25th and 75th percentiles, and triangles are minimum and maximum values. LQR approximation accuracy is shown as a dashed line.

$V$-NN and $V$-*QRnet* are trained as described in [8]. For the value gradient networks we do not use the value gradient loss term (23) since for this problem it did not improve results. To be consistent, all NNs have $L = 5$ hidden layers with $w = 32$ neurons each and $\tanh(\cdot)$ nonlinearities. Optimization of (22) is carried out with L-BFGS [37], which stops when the relative change in the loss is sufficiently small. All models are trained on an NVIDIA RTX 2080Ti GPU.

Figure 3 shows training times and test accuracies of the NNs. We see that models which approximate the value gradient, especially $\lambda_{\text{Jac}}$-*QRnet* and $\lambda_{\text{mat}}$-*QRnet*, take the longest to train because of the large number of NN parameters. Despite this, the training time is still very reasonable at under eight minutes. We also find that the new architectures have similar test accuracy statistics to the standard NNs, confirming that they can learn complicated nonlinear functions as suggested by Theorems 1 and 2. For this problem there is no clear performance distinction between the "Jacobian" and "matrix" architectures or between the $\lambda$ and $u$ models.

### B. LOCAL STABILITY VERIFICATION

As a first step we assess the local stability of each NN-controlled system. Let $\bar{\mathbf{x}} \in \mathbb{R}^n$ be an equilibrium of the closed loop system, $\dot{\mathbf{x}} = \mathbf{f}(\mathbf{x}, \widehat{\mathbf{u}}(\mathbf{x}))$, i.e. $\mathbf{f}(\bar{\mathbf{x}}, \widehat{\mathbf{u}}(\bar{\mathbf{x}})) = \mathbf{0}$. Note that the NN controllers introduced in Section III always have $\mathbf{f}(\mathbf{x}_f, \widehat{\mathbf{u}}(\mathbf{x}_f)) = \mathbf{0}$, but $\mathbf{x}_f$ may not be a closed loop equilibrium for $V$-NN, $\lambda$-NN, $u$-NN, and $V$-*QRnet*. Thus for these controllers we use a root-finding algorithm to locate a closed loop equilibrium $\bar{\mathbf{x}}$ near $\mathbf{x}_f$. The dynamics near $\bar{\mathbf{x}}$ can be approximated by $\dot{\mathbf{x}} \approx \mathbf{A}_{\text{cl}}(\mathbf{x} - \bar{\mathbf{x}})$, where

$$\mathbf{A}_{\text{cl}} := \frac{\partial \mathbf{f}}{\partial \mathbf{x}}(\bar{\mathbf{x}}, \widehat{\mathbf{u}}(\bar{\mathbf{x}})) + \left[\frac{\partial \mathbf{f}}{\partial \mathbf{u}}(\bar{\mathbf{x}}, \widehat{\mathbf{u}}(\bar{\mathbf{x}}))\right]\left[\frac{\partial \widehat{\mathbf{u}}}{\partial \mathbf{x}}(\bar{\mathbf{x}})\right] \quad (26)$$

is the closed loop Jacobian. Therefore, after synthesizing a feedback controller we can easily verify local stability by seeing if $\mathbf{A}_{\text{cl}}$ is Hurwitz. Furthermore as noted in [26], one benefit of using an NN controller with differentiable activation functions is that the closed loop dynamics are locally $\mathcal{C}^1$. This allows one to use tools from linear systems theory to characterize the local stability of $\bar{\mathbf{x}}$.

Figure 4 shows the real part of the most positive eigenvalue of $\mathbf{A}_{\text{cl}}$ for each NN. We find that standard NNs must be trained to a high level of test accuracy before they are even locally stable, which necessitates a large data set and long training time. On the other hand, $\lambda$-*QRnet*. $u$-*QRnet*, and the new "Jacobian" and "matrix" *QRnet*s all yield LAS even when trained on small data sets. Recall that Proposition 1 guarantees this for the new architectures.

### C. MONTE CARLO STABILITY ANALYSIS

Here and in Section D we conduct Monte Carlo (MC) closed loop simulations. We randomly select $N_{\text{MC}} = 100$ initial conditions $\mathbf{x}_0^{(i)}$, $i = 1, \ldots, N_{\text{MC}}$ with norm $\left\|\mathbf{x}_0^{(i)} - \mathbf{x}_f\right\| = 1.2 \approx \max_{\mathbf{x}^{(j)} \in \mathcal{D}_{\text{train}}}\left\|\mathbf{x}^{(j)} - \mathbf{x}_f\right\|$, placing them at the edge of the training domain where the NNs may be less accurate and the system harder to control. We stop each simulation when the system reaches a steady state or exceeds a large final time. We call the largest observed final state,

$$\max_{\mathbf{x}_0^{(i)}} \left\|\mathbf{x}\left(t_f; \mathbf{x}_0^{(i)}\right) - \mathbf{x}_f\right\|,$$

the *worst-case failure*. If this is sufficiently small then the closed loop nonlinear system is likely to be semi-globally stable.

Figure 5 shows the worst-case failures for each controller. We find that only the most accurate standard NNs stabilize the origin, whereas *all* controllers from Section III stabilize all the MC trajectories. These empirical results suggest that the proposed architectures not only guarantee LAS, but also make the control design process more reliable, consistently yielding a stabilizing control law even with small data sets and short training times.

### D. MONTE CARLO OPTIMALITY ANALYSIS

In this paper we are interested in both stability and optimality. Optimality of a given controller $\widehat{\mathbf{u}}(\cdot)$ can be quantified by the accumulated cost $J\left[\widehat{\mathbf{u}}(\cdot); \mathbf{x}_0^{(i)}\right]$ compared to the optimal cost $V\left(\mathbf{x}_0^{(i)}\right)$, across all MC simulations $i = 1, \ldots, N_{\text{MC}}$. Figure 6 shows the results of this analysis for the same set of MC simulations conducted in Section C. Among the stabilizing NN controllers there is a clear correlation between higher test accuracy and better performance. All the stabilizing NN controllers follow this trend and perform better than LQR alone. It follows that the proposed architectures improve stability without limiting optimality.

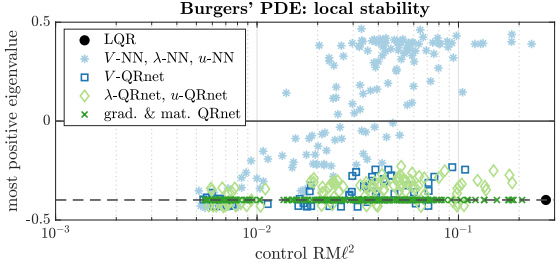

Figure 4: Real part of most positive eigenvalue of the closed loop Jacobian $\mathbf{A}_{\mathrm{cl}}$ at or near $\mathbf{x}_f$. Each marker represents a single model.

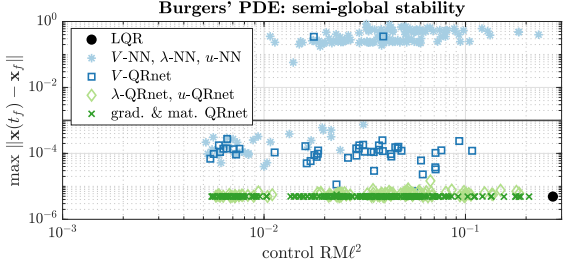

Figure 5: Worst-case norm of final state over $N_{\mathrm{MC}} = 100$ simulations.

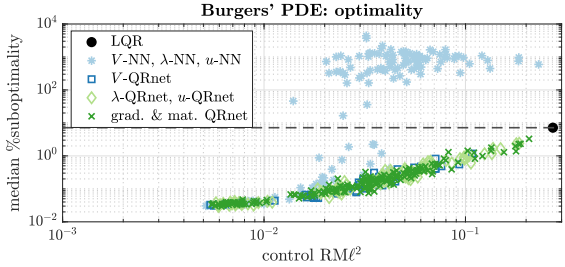

Figure 6: Median percent cost more than optimal cost over $N_{\mathrm{MC}} = 100$ simulations.

## VI. APPLICATION EXAMPLE: FIXED-WING UAV

In this section we illustrate how the proposed control architectures can be used with supervised learning to design an optimal feedback controller a fixed-wing 6DoF UAV. The controller is designed for stabilization from a wide range of flight conditions, as well as tracking for arbitrary altitude and course commands - a challenging nonlinear OCP.

### A. FIXED-WING UAV DYNAMICS

The dynamic model we use is based on the one presented in [38], [39]. We review it here to orient the reader and point out several small differences.

The position of the UAV is described in inertial north-east-down coordinates, $\mathbf{p} := \begin{pmatrix} p_n, & p_e, & p_d \end{pmatrix}^T$. The velocities in the body $x$, $y$, and $z$ directions are denoted as $\mathbf{V} := \begin{pmatrix} u, & v, & w \end{pmatrix}^T$. The attitude of the UAV, i.e. its rotation from inertial to body frames, is described using quaternions $\mathbf{q} := \begin{pmatrix} q_0, & \bar{\mathbf{q}}^T \end{pmatrix}^T$, where $q_0$ is the scalar quaternion and $\bar{\mathbf{q}} := \begin{pmatrix} q_1, & q_2, & q_3 \end{pmatrix}^T$ is the vector quaternion. The angular velocity of UAV in the body frame is written as $\boldsymbol{\omega} := \begin{pmatrix} p, & q, & r \end{pmatrix}^T$. The full state is then

$$\mathbf{x} := \begin{pmatrix} \mathbf{p}^T, & \mathbf{V}^T, & \mathbf{q}^T, & \boldsymbol{\omega}^T \end{pmatrix}^T \in \mathbb{R}^{13}. \quad (27)$$

The UAV is controlled with a throttle $\delta_t \in [0, 1]$, ailerons $\delta_a \in [-\delta_a^+, \delta_a^+]$, elevator $\delta_e \in [-\delta_e^+, \delta_e^+]$, and rudder $\delta_r \in [-\delta_r^+, \delta_r^+]$. Thus

$$\mathbf{u} := \begin{pmatrix} \delta_t, & \delta_a, & \delta_e, & \delta_r \end{pmatrix}^T \in \mathbb{U} \subset \mathbb{R}^4. \quad (28)$$

Modeling the UAV as a rigid body we obtain the dynamic equations [38]

$$\dot{\mathbf{x}} = \begin{pmatrix} \dot{\mathbf{p}} \\ \dot{\mathbf{V}} \\ \dot{\mathbf{q}} \\ \dot{\boldsymbol{\omega}} \end{pmatrix} = \begin{pmatrix} \mathcal{R}_{\mathbf{q}}^{-1}(\mathbf{V}) \\ -\boldsymbol{\omega} \times \mathbf{V} + \frac{1}{m}\mathbf{F} \\ \frac{1}{2}\boldsymbol{\omega}\mathbf{q} \\ \mathbf{J}^{-1}\left[-\boldsymbol{\omega} \times (\mathbf{J}\boldsymbol{\omega}) + \mathbf{M}\right] \end{pmatrix}. \quad (29)$$

Here $\mathcal{R}_{\mathbf{q}} : \mathbb{R}^3 \to \mathbb{R}^3$ is the rotation (computed using the attitude $\mathbf{q}$) from inertial to body frame, and $\mathcal{R}_{\mathbf{q}}^{-1}(\cdot)$ is the inverse rotation from body to inertial frame. Next, $m$ is the mass of the UAV, $\mathbf{J} \in \mathbb{R}^{3 \times 3}$ is the UAV's inertia matrix, and

$$\boldsymbol{\omega} := \begin{pmatrix} 0 & -p & -q & -r \\ p & 0 & r & -q \\ q & -r & 0 & p \\ r & q & -p & 0 \end{pmatrix}.$$

Finally $\mathbf{F} = \mathbf{F}(\mathbf{x}, \mathbf{u})$ and $\mathbf{M} = \mathbf{M}(\mathbf{x}, \mathbf{u})$ are the external forces and moments acting on the vehicle expressed in the body frame. These arise as a result of gravity, aerodynamics, and control inputs. In the following presentation we ignore the effects of wind for simplicity.

The first force is gravity, which can be expressed in the body frame as $\mathbf{F}_{\mathrm{gravity}} = \mathcal{R}_{\mathbf{q}} \left[\begin{pmatrix} 0, & 0, & mg \end{pmatrix}^T\right]$, where $g$ is the gravitational constant. Next we employ a linear propeller model based on [38]:

$$\mathbf{F}_{\mathrm{prop}} = \frac{1}{2}\rho\pi R_{\mathrm{prop}}^2 C_{\mathrm{prop}} \begin{pmatrix} k_{\mathrm{motor}}^2 \delta_t - \|\mathbf{V}\|^2 \\ 0 \\ 0 \end{pmatrix}. \quad (30)$$

Here $\rho$ is the air density, $R_{\mathrm{prop}}$ is the propeller blade length, and $k_{\mathrm{motor}}$ and $C_{\mathrm{prop}}$ parameterize thrust efficiency.

Finally, the aerodynamic forces $\mathbf{F}_{\mathrm{aero}} = \begin{pmatrix} F_x, & F_y, & F_z \end{pmatrix}$ and moments $\mathbf{M} = \begin{pmatrix} M_\ell, & M_m, & M_n \end{pmatrix}$ are in general complicated nonlinear relationships that must be modeled from experimental data. In this work we use the basic models from [38], with slight nonlinear modifications to the drag and pitching moment models to improve their post-stall realism. The longitudinal forces are modeled as

$$\begin{pmatrix} F_x \\ F_z \end{pmatrix} = \begin{pmatrix} \cos\alpha & -\sin\alpha \\ \sin\alpha & \cos\alpha \end{pmatrix} \begin{pmatrix} -F_D \\ -F_L \end{pmatrix},$$

where $\alpha = \tan^{-1}(w/u)$ is the angle of attack and

$$F_L = \frac{1}{2}\rho\|\mathbf{V}\|^2 S \left[C_L(\alpha) + \frac{cC_{L_q}}{2\|\mathbf{V}\|}q + C_{L_{\delta_e}}\delta_e\right],$$

$$F_D = \frac{1}{2}\rho\|\mathbf{V}\|^2 S \left[ C_D(\alpha) + \frac{cC_{D_q}}{2\|\mathbf{V}\|}q + C_{D_{\delta_e}}\delta_e \right],$$

are the lift and drag forces, respectively. Here $S$ is the wing area and $C_{L_q}$, $C_{L_{\delta_e}}$, $C_{D_q}$, $C_{D_{\delta_e}}$ are modeling parameters. As in [38] we model

$$C_L(\alpha) = [1-\sigma_b(\alpha)][C_{L_0} + C_{L_\alpha}\alpha] \\ + \sigma_b(\alpha) \cdot 2\mathrm{sign}(\alpha)\sin^2\alpha\cos\alpha, \quad (31)$$

where $C_{L_0}$ and $C_{L_\alpha}$ are modeling parameters and $\sigma_b(\alpha)$ is a smooth blending function which is $\sigma_b(\alpha) \approx 0$ for $|\alpha| < \alpha_{\mathrm{stall}}$ and $\sigma_b(\alpha) \approx 1$ for $|\alpha| > \alpha_{\mathrm{stall}}$, with $\alpha_{\mathrm{stall}}$ being the stall angle of attack. See [38] for details. For the drag model we use a blend of a quadratic and post-stall flat plate model [40]:

$$C_D(\alpha) = [1-\sigma_b(\alpha)]\left[ C_{D_0} + \frac{(C_{L_0}+C_{L_\alpha}\alpha)^2}{\pi e b^2 /S}\right] \\ + \sigma_b(\alpha) \cdot 2\sin^2\alpha, \quad (32)$$

where $C_{D_0}$ is the parasitic drag, $b$ is the wingspan, and $e$ is another modeling parameter. We similarly modify the pitching moment model from [38] to be nonlinear in $\alpha$. Let

$$M_m = \frac{1}{2}\rho\|\mathbf{V}\|^2 Sc\left[C_m(\alpha) + \frac{C_{m_q}c}{2\|\mathbf{V}\|}q + C_{m_{\delta_e}}\delta_e\right], \quad (33)$$

with

$$C_m(\alpha) = [1-\sigma_b(\alpha)]\tanh(C_{m_0} + C_{m_\alpha}\alpha) \\ + \sigma_b(\alpha) \cdot C_{m_\infty}\sin(-\alpha), \quad (34)$$

and where $C_{m_q}$, $C_{m_{\delta_e}}$, $C_{m_0}$, $C_{m_\alpha}$, and $C_{m_\infty}$ are modeling parameters. The remaining lateral aerodynamics, $F_y$, $M_\ell$, and $M_n$, are functions of $\|\mathbf{V}\|$, $p$, $r$, $\delta_a$, $\delta_r$, and the sideslip $\beta = \sin^{-1}(v/\|\mathbf{V}\|)$. These models are the same as in [38] and are omitted for brevity. The values of the constants used in this problem are taken from [39], with the exception of $C_{\mathrm{prop}} = 0.45$, $k_{\mathrm{motor}} = 32$, $\alpha_{\mathrm{stall}} = 20°$, and $C_{m_\infty} = 0.8$. Note that $\alpha_{\mathrm{stall}} = 20°$ is lower than in [38], [39], making the model more realistic and challenging to control.

### B. OPTIMAL CONTROL PROBLEM FORMULATION

We aim to design a feedback controller to stabilize the UAV and track any desired altitude $h_f = -p_{d,f}$ and course angle $\chi_f = \tan^{-1}(\dot{p}_e/\dot{p}_n)$. Let $\mathbf{x}_f, \mathbf{u}_f$ be the pair of *trim* states and controls computed for a desired airspeed $\|\mathbf{V}_f\|$. The UAV is in trim if $\mathbf{f}(\mathbf{x}_f, \mathbf{u}_f) = \mathbf{0}$, except for $\dot{p}_n$ and $\dot{p}_e$. Note that the dynamics are invariant to $\mathbf{p}$, so we can choose any arbitrary trim altitude. The dynamics (excepting $\dot{p}_n$ and $\dot{p}_e$) are also invariant to rotations of the inertial reference frame about the inertial $z$ axis, which allows us to use the same trim attitude $\mathbf{q}_f$ to express any desired yaw angle $\psi_f$. When the vehicle is in trim (and in the absence of wind) the yaw angle $\psi$ is equal to the course angle $\chi$, and thus this formulation allows arbitrary course tracking.

A suitable running cost for this OCP is

$$\mathcal{L}(\mathbf{x},\mathbf{u}) = Q_h\left[h_{\mathrm{ceil}}\tanh\left(\frac{p_d - p_{d,f}}{h_{\mathrm{ceil}}}\right)\right]^2 \\ + (\mathbf{V}-\mathbf{V}_f)^T \mathbf{Q_V}(\mathbf{V}-\mathbf{V}_f)$$

$$+ (\bar{\mathbf{q}}-\bar{\mathbf{q}}_f)^T \mathbf{Q_q}(\bar{\mathbf{q}}-\bar{\mathbf{q}}_f)^T \\ + (\boldsymbol{\omega}-\boldsymbol{\omega}_f)^T \mathbf{Q_\omega}(\boldsymbol{\omega}-\boldsymbol{\omega}_f) \\ + (\mathbf{u}-\mathbf{u}_f)^T \mathbf{R}(\mathbf{u}-\mathbf{u}_f), \quad (35)$$

where $Q_h, h_{\mathrm{ceil}} > 0$, $\mathbf{Q_V}, \mathbf{Q_q}, \mathbf{Q_\omega} \in \mathbb{R}^{3\times 3}$ are positive definite, and $\mathbf{R} \in \mathbb{R}^{4\times 4}$ is positive definite. Notice that the altitude cost is locally quadratic but saturates for $|p_d - p_{d,f}| \geq h_{\mathrm{ceil}}$, preventing extreme maneuvers when the commanded altitude changes.

We set the desired airspeed at $\|\mathbf{V}_f\| = 20$ [m/s] and use the following cost function parameters:

$$\begin{cases} h_{\mathrm{ceil}} = 50 \text{ [m]}, \quad Q_h = 1/h_{\mathrm{ceil}}^2, \\ \mathbf{Q_V} = \mathrm{diag}\left(10/\|\mathbf{V}_f\|^2, \ 1, \ 1\right), \\ \mathbf{Q_q} = 5\mathbf{I}_{3\times 3}, \\ \mathbf{Q_\omega} = \frac{1}{(30°)^2}\mathbf{I}_{3\times 3}, \\ \mathbf{R} = \mathrm{diag}\left(0.1, \ 0.1/(\delta_a^+)^2, \ 1/(\delta_e^+)^2, \ 1/(\delta_r^+)^2\right). \end{cases}$$

Initial conditions are uniformly sampled from the following domain to elicit a wide range of nonlinear dynamics:

$$\mathbb{X}_0 = \left\{ \begin{array}{l} p_{d_0} \in [-3h_{\mathrm{ceil}}, 3h_{\mathrm{ceil}}], \\ \mathbf{V}_0 \in [\mathbf{V}_f - 5 \text{ [m/s]}, \mathbf{V}_f + 5 \text{ [m/s]}], \\ \psi_0, \phi_0 \in [-180°, 180°], \quad \theta_0 \in [-90°, 90°], \\ \boldsymbol{\omega}_0 \in [-30 \text{ [deg/s]}, 30 \text{ [deg/s]}]. \end{array} \right\}$$

Here $\psi_0, \theta_0, \phi_0$ denote the inital yaw, pitch, and roll angles, which are converted to the initial quaternion $\mathbf{q}_0$. Recall that we can set $p_{d_f} = 0$ and $\psi_f = 0$ without loss of generality, and thus the initial condition determines the initial altitude and course errors.

For this high-dimensional and highly nonlinear OCP we found that indirect methods were unreliable for generating data. For this reason we generate data with a direct method: Radau pseudospectral method [33]. To the best of our knowledge this is the first case of pseudospectral methods being used for supervised learning. To obtain good quality open loop OCP data we use a large number of Radau collocation points and set stringent tolerances for the nonlinear programming solver. Note that direct methods typically provide optimal state and control pairs, but obtaining costate entails extra computational effort. For this reason in this section we only show results for $u$-NN, $u$-QRnet, $u_{\mathrm{Jac}}$-QRnet, and $u_{\mathrm{mat}}$-QRnet, which directly approximate the optimal control and do not need costate data.

### C. LEARNING RESULTS

For this problem we generate training data sets with 32, 64, 128, and 256 trajectories each. For each data set we train each type of NN controller with different weight initializations. As before we conduct ten trials for each data set size. We evaluate the RM$\ell^2$ error (24) on an independent data set with 100 trajectories. As in Section V all NNs have $L = 5$ hidden layers with $w = 32$ neurons each and $\tanh(\cdot)$ nonlinearities. Because these data sets are too large for full-batch optimization we use Adam [41] with a learning rate of $10^{-3}$, batch sizes of 256 data points, and 1500 epochs.

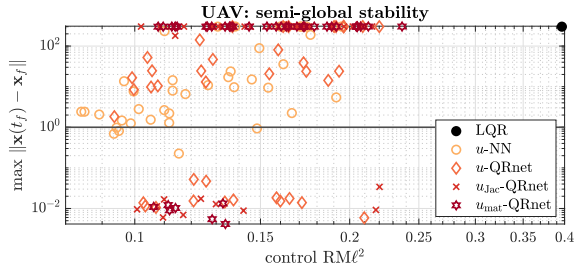

Figure 7: Worst-case norm of final state over $N_{\mathrm{MC}} = 100$ simulations. The vertical axis is limited since we stop simulations if the altitude $h - h_f$ goes outside of $\pm 300$ [m].

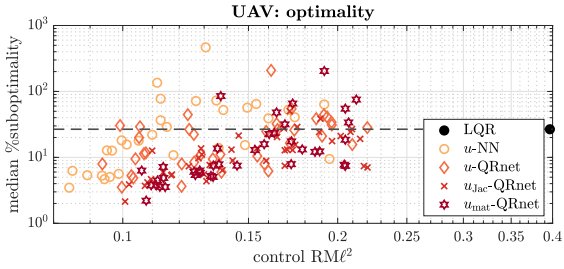

Figure 8: Median percent cost more than optimal cost over $N_{\mathrm{MC}} = 100$ simulations.

Similar to the results shown in Figure 4, again we find that well-trained $u$-NNs may fail to even locally stabilize the system. Furthermore, closed loop equilibria under $u$-NN control are often far from $\mathbf{x}_f$. In the physical system this corresponds to steady state altitude, course, and attitude errors, even when said equilibrium is stable (see Figure 1).

Figure 7 shows the worst case norm for a set of $N_{\mathrm{MC}} = 100$ closed loop simulations. These simulations demonstrate how challenging the UAV is to control over this large spatial domain. First we notice that LQR is *not* globally stabilizing for this OCP. Next we observe that most standard $u$-NNs, even the well-trained ones, do not stabilize $\mathbf{x}_f$. $u$-QRnet, $u_{\mathrm{Jac}}$-QRnet, and $u_{\mathrm{mat}}$-QRnet also have some difficulty with semi-global stabilization, though they clearly do better than $u$-NN. Note that these controllers are able to stabilize trajectories LQR fails to stabilize - even though they are built on top of LQR.

Finally Figure 8 shows the average performance of each controller in terms of minimizing the cost functional $J[\mathbf{u}(\cdot)]$. We again see that most NN controllers perform better than LQR on average, indicating they they do learn the optimal policy reasonably well. We also see that the standard $u$-NNs have slightly higher test accuracy, suggesting that for this OCP the training loss (22) converges faster than the modified architecture (i.e. requires fewer gradient descent steps). Despite this, we can see that $u$-QRnet, $u_{\mathrm{Jac}}$-QRnet, and $u_{\mathrm{mat}}$-QRnet perform just as good or better in terms of closed loop stability and optimality. We expect that all methods will improve with larger data sets, more training epochs, and hyperparameter tuning.

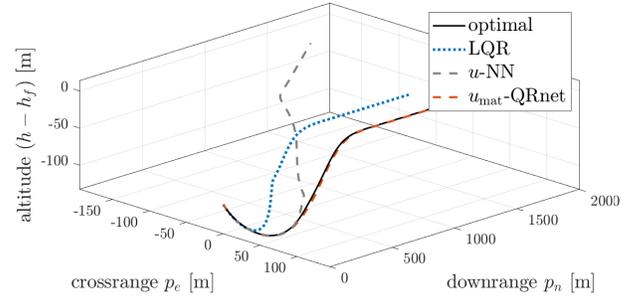

Figure 9: Simulated trajectory of the fixed wing UAV (29) with $u$-NN, $u_{\mathrm{mat}}$-QRnet, and LQR controllers, compared to optimal trajectory.

### D. EXAMPLE CLOSED LOOP SIMULATION

We conclude this section with an illustrative example of an NN-in-the-loop simulation. We conduct the simulation for the same initial condition as in Figure 1, this time with a $u_{\mathrm{mat}}$-QRnet controller trained on 256 trajectories. A view of the closed loop trajectory is presented in Figure 9, and detailed time series of system states and feedback controls are shown in Figure 10. Notice that the UAV begins off course and pitched down with large negative pitch rate. For this initial condition, LQR is 36.44% suboptimal while the $u_{\mathrm{mat}}$-QRnet is 0.95% suboptimal. This simulation highlights the potential benefit of using NN optimal feedback controllers to achieve good performance in nonlinear systems.

## VII. CONCLUSION

In this paper we have shown that NN feedback controllers can frequently fail to stabilize a system, even when they are trained to a high degree of accuracy. This occurs frequently enough that it cannot be ignored. One strategy to make NN feedback controllers more viable is through the use of specialized NN architectures. To this end we propose four new model architectures which guarantee (at least) LAS while retaining the approximation capacity necessary to learn the full nonlinear optimal control and provide nonlinear stability on semi-global domains. A summary of the control architectures discussed in this paper is given in Table 1.

In Section V we evaluated the proposed architectures through a series of practical closed loop stability and optimality tests, demonstrating their advantages over standard NNs. Finally in Section VI we illustrated how the proposed architectures might be used with supervised learning to design optimal feedback controllers for challenging, practical systems.

For problems where the dimension is not too large, the value gradient approximators, $\lambda_{\mathrm{Jac}}$-QRnet and $\lambda_{\mathrm{mat}}$-QRnet, can sometimes perform better than the control approximators, $u_{\mathrm{Jac}}$-QRnet and $u_{\mathrm{mat}}$-QRnet. This is because they encode additional physical structure and can learn from costate data in addition to control data. On the other hand, the control approximators are generally much faster to train, and they can be implemented even when it is not possible to solve (6) for the optimal control, and when it is difficult to generate accurate costate data.

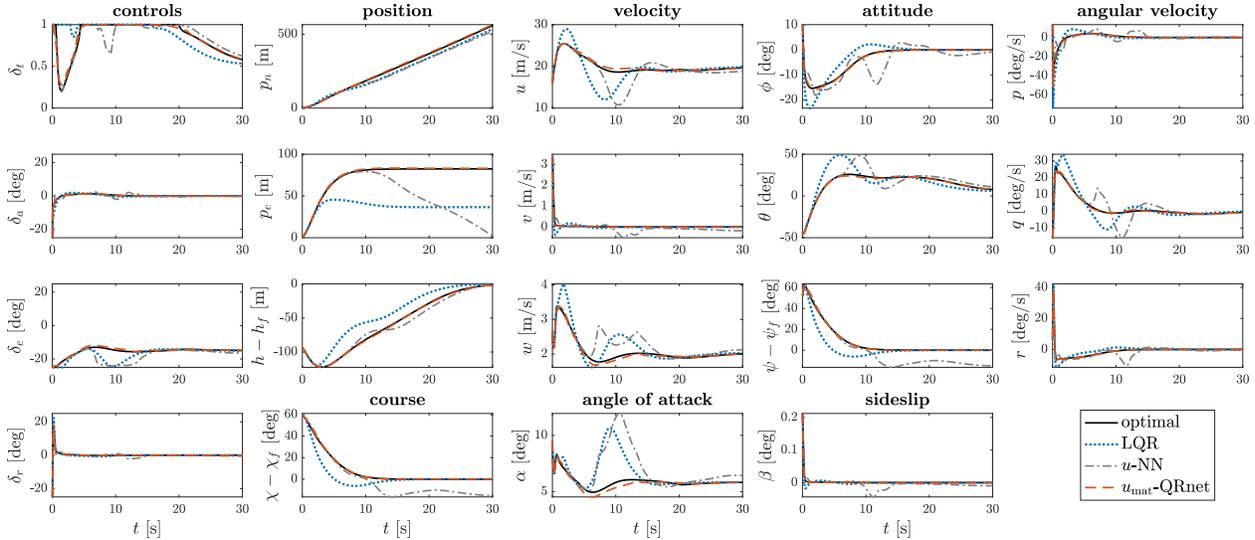

Figure 10: Closed loop simulations with $u$-NN, $u_{\text{mat}}$-*QRnet*, and LQR controllers compared to optimal trajectory and controls. Attitude is given in Euler angles $\phi, \theta, \psi$ (roll, pitch, yaw).

| architecture | LAS guarantee | number of NN parameters ($p$) |
|---|---|---|
| $V$-NN* | - | $2w + Lw^2$ |
| $V$-*QRnet** [8] | - | $2w + Lw^2$ |
| $\lambda$-NN* | - | $2wn + Lw^2$ |
| $\lambda$-*QRnet** (11) | - | $2wn + Lw^2$ |
| $\lambda_{\text{Jac}}$-*QRnet** (17) | ✓ | $2wn + Lw^2$ |
| $\lambda_{\text{mat}}$-*QRnet** (19) | ✓ | $wn^2 + wn + Lw^2$ |
| $u$-NN | - | $wm + wn + Lw^2$ |
| $u$-*QRnet* (13) | - | $wm + wn + Lw^2$ |
| $u_{\text{Jac}}$-*QRnet* (18) | ✓ | $wm + wn + Lw^2$ |
| $u_{\text{mat}}$-*QRnet* (20) | ✓ | $wmn + wn + Lw^2$ |

Table 1: Summary of NN control architectures discussed in this paper. $L$ denotes the number of layers and $w$ is their width. *Need to solve (6) for optimal control in terms of state and costate.

To complement the NN architectures presented in this paper, in future work we intend to develop mathematical tools to explain the behavior of NN feedback controllers. In particular, we would like to better understand what causes seemingly-accurate NN models to fail at stabilizing a system, as well as why and to what extent the novel NN architectures improve semi-global system stability. Such theoretical advances will be necessary if supervised learning is to become a reliable and commonly accepted control design method.

## APPENDIX: PROOFS OF APPROXIMATION CAPACITY

To prove the approximation capacity theorems, we first specialize the classic Stone-Weierstrass theorem (stated below for reference) to approximation of locally $\mathcal{C}^1$ functions. This result is given as Corollary 1. The proofs of Theorems 1 and 2 in Sections A and B, respectively, subsequently apply classical NN approximation theory [28] to the approximating function from Corollary 1 to obtain the desired result.

By $\mathcal{C}(\mathbb{X})$ and $\mathcal{C}(\mathbb{X}; \mathbb{R}^d)$ we denote the spaces of continuous functions on $\mathbb{X}$ taking values in $\mathbb{R}$ and $\mathbb{R}^d$, respectively. These function spaces are *algebras*; a set of functions $\mathcal{A}$ an algebra if it is closed under (element-wise) addition, multiplication, and scalar multiplication. A subalgebra of $\mathcal{A}$ is a subset of $\mathcal{A}$ which is also an algebra.

**Theorem 3 (Stone-Weierstrass [29]):**
Suppose that $\mathcal{A}$ is a subalgebra of $\mathcal{C}(\mathbb{X})$ which separates points[4] and does not vanish[5] anywhere in $\mathbb{X}$. Then for all $f \in \mathcal{C}(\mathbb{X})$ and all $\varepsilon > 0$ there exists $g \in \mathcal{A}$ satisfying $\max_{\mathbf{x} \in \mathbb{X}} |f(\mathbf{x}) - g(\mathbf{x})| < \varepsilon$.

**Corollary 1 (Approximation of locally $\mathcal{C}^1$ functions):**
Suppose $\mathbf{f} \in \mathcal{C}(\mathbb{X}; \mathbb{R}^d)$, $f(\mathbf{0}) = \mathbf{0}$, and $\mathbf{f}(\cdot)$ is $\mathcal{C}^1$ in a neighborhood of $\mathbf{0}$. Then for all $\varepsilon > 0$, there exists a function $\mathbf{g} \in \mathcal{C}^1(\mathbb{X}; \mathbb{R}^d)$ satisfying $\mathbf{g}(\mathbf{0}) = \mathbf{0}$, $[\partial \mathbf{g}/\partial \mathbf{x}](\mathbf{0}) = \mathbf{0}$, and $\max_{\mathbf{x} \in \mathbb{X}} \left\| \mathbf{f}(\mathbf{x}) - \left[\frac{\partial \mathbf{f}}{\partial \mathbf{x}}(\mathbf{0})\right] \mathbf{x} - \mathbf{g}(\mathbf{x}) \right\|_1 < \varepsilon$.

*Proof:* Consider the set of functions $\mathcal{A} \subset \mathcal{C}^1(\mathbb{X})$ which have $[\partial g/\partial \mathbf{x}](\mathbf{0}) = \mathbf{0}$. We claim that $\mathcal{A}$ is an algebra which vanishes nowhere and separates points. It is easy to verify that $\mathcal{A}$ is closed under addition, multiplication, and scalar multiplication; hence $\mathcal{A}$ is an algebra. $\mathcal{A}$ also contains the constant functions, so it vanishes nowhere. Lastly, to see that $\mathcal{A}$ separates points, note that for any $\mathbf{x} \neq \mathbf{y}$ without loss of generality $x_1 \neq y_1$. Then take $g(\mathbf{x}) = x_1^3 \neq y_1^3 = g(\mathbf{y})$, since $x^3$ is one-to-one.

---

[4] An algebra $\mathcal{A}$ separates points if for all $\mathbf{x}, \mathbf{y} \in \mathbb{X}$, $\mathbf{x} \neq \mathbf{y}$, there exists $g \in \mathcal{A}$ such that $g(\mathbf{x}) \neq g(\mathbf{y})$.

[5] A set of functions $\mathcal{A}$ vanishes at $\mathbf{x}_1 \in \mathbb{X}$ if $\mathbf{f}(\mathbf{x}_1) = \mathbf{0}$ for all $\mathbf{f} \in \mathcal{A}$.

Now write $\mathbf{f}(\mathbf{x}) = \left(f_1(\mathbf{x}), \ldots, f_d(\mathbf{x})\right)^T$. For each $i = 1, \ldots, d$, by Theorem 3 there exists $h_i \in \mathcal{A}$ satisfying

$$\max_{\mathbf{x} \in \mathbb{X}} \left| f_i(\mathbf{x}) - \left[\tfrac{\partial f_i}{\partial \mathbf{x}}(\mathbf{0})\right] \mathbf{x} - h_i(\mathbf{x}) \right| < \tfrac{\varepsilon}{2d}.$$

Since $f_i(\mathbf{0}) = 0$ by assumption, this implies

$$|h_i(\mathbf{0})| = \left| f_i(\mathbf{0}) - \left[\tfrac{\partial f_i}{\partial \mathbf{x}}(\mathbf{0})\right] \mathbf{0} - h_i(\mathbf{0}) \right| < \tfrac{\varepsilon}{2d}.$$

Defining $g_i(\mathbf{x}) := h_i(\mathbf{x}) - h_i(\mathbf{0})$ we get $g_i(\mathbf{0}) = 0$ and

$$\max_{\mathbf{x} \in \mathbb{X}} \left| f_i(\mathbf{x}) - \left[\tfrac{\partial f_i}{\partial \mathbf{x}}(\mathbf{0})\right] \mathbf{x} - g_i(\mathbf{x}) \right|$$
$$\leq \max_{\mathbf{x} \in \mathbb{X}} \left| f_i(\mathbf{x}) - \left[\tfrac{\partial f_i}{\partial \mathbf{x}}(\mathbf{0})\right] \mathbf{x} - h_i(\mathbf{x}) \right| + |h_i(\mathbf{0})|$$
$$< \tfrac{\varepsilon}{2d} + \tfrac{\varepsilon}{2d} = \tfrac{\varepsilon}{d}.$$

Because $\mathcal{A}$ is an algebra we must also have $g_i \in \mathcal{A}$ and hence $[\partial g_i / \partial \mathbf{x}](\mathbf{0}) = \mathbf{0}$. Thus setting $\mathbf{g}(\mathbf{x}) = \left(g_1(\mathbf{x}), \ldots, g_d(\mathbf{x})\right)^T$ yields the desired function. ∎

### A. PROOF OF THEOREM 1

Since $\mathbb{X}$ is bounded, there is some $B > 0$ for which $\max_{\mathbf{x} \in \mathbb{X}} \|\mathbf{x}\|_1 \leq B$. For any $\varepsilon > 0$, let $\varepsilon^* = \min\{\varepsilon, \varepsilon/B\}$. From Corollary 1 we can find some $\mathbf{g} \in \mathcal{C}^1\left(\mathbb{X}; \mathbb{R}^d\right)$ satisfying $\mathbf{g}(\mathbf{0}) = \mathbf{0}$, $[\partial \mathbf{g}/\partial \mathbf{x}](\mathbf{0}) = \mathbf{0}$, and

$$\max_{\mathbf{x} \in \mathbb{X}} \left\| \mathbf{f}(\mathbf{x}) - \left[\tfrac{\partial \mathbf{f}}{\partial \mathbf{x}}(\mathbf{0})\right] \mathbf{x} - \mathbf{g}(\mathbf{x}) \right\|_1 < \tfrac{\varepsilon^*}{2}.$$

Also by the universal approximation theorem [28] there exists an NN, $\mathcal{N} \in \mathcal{C}^1\left(\mathbb{X}; \mathbb{R}^d\right)$, which approximates $\mathbf{g}(\cdot)$ and its derivative to arbitrary accuracy, say

$$\max \left\{ \begin{array}{c} \max_{\mathbf{x} \in \mathbb{X}} \|\mathbf{g}(\mathbf{x}) - \mathcal{N}(\mathbf{x})\|_1, \\ \max_{\mathbf{x} \in \mathbb{X}} \left\| \tfrac{\partial \mathbf{g}}{\partial \mathbf{x}}(\mathbf{x}) - \tfrac{\partial \mathcal{N}}{\partial \mathbf{x}}(\mathbf{x}) \right\|_{1,1} \end{array} \right\} < \tfrac{\varepsilon^*}{6}.$$

Here we define the matrix norm $\|\mathbf{A}\|_{1,1} := \|\text{vec}(\mathbf{A})\|_1$ for a matrix $\mathbf{A} \in \mathbb{R}^{d \times n}$ and its vectorization, $\text{vec}(\mathbf{A}) \in \mathbb{R}^{mn}$. Notice that

$$\|\mathcal{N}(\mathbf{0})\|_1 = \|\mathbf{0} - \mathcal{N}(\mathbf{0})\|_1 = \|\mathbf{g}(\mathbf{0}) - \mathcal{N}(\mathbf{0})\|_1 < \tfrac{\varepsilon^*}{6}.$$

Consequently, for all $\mathbf{x} \in \mathbb{X}$ we get

$$\|\mathbf{g}(\mathbf{x}) - \mathcal{N}(\mathbf{x}) - \mathcal{N}(\mathbf{0})\|_1$$
$$\leq \|\mathbf{g}(\mathbf{x}) - \mathcal{N}(\mathbf{x})\|_1 + \|\mathcal{N}(\mathbf{0})\|_1 < \tfrac{\varepsilon^*}{3}.$$

Similarly,

$$\left\| \tfrac{\partial \mathbf{g}}{\partial \mathbf{x}}(\mathbf{0}) - \tfrac{\partial \mathcal{N}}{\partial \mathbf{x}}(\mathbf{0}) \right\|_{1,1} = \left\| \tfrac{\partial \mathcal{N}}{\partial \mathbf{x}}(\mathbf{0}) \right\|_{1,1} < \tfrac{\varepsilon^*}{6},$$

which implies

$$\left\| \left[\tfrac{\partial \mathcal{N}}{\partial \mathbf{x}}(\mathbf{0})\right] \mathbf{x} \right\|_1 \leq \left\| \tfrac{\partial \mathcal{N}}{\partial \mathbf{x}}(\mathbf{0}) \right\|_{1,1} \|\mathbf{x}\|_1 < \tfrac{\varepsilon^*}{6} B,$$

for all $\mathbf{x} \in \mathbb{X}$. Putting this all together we obtain

$$\left\| \mathbf{f}(\mathbf{x}) - \left(\left[\tfrac{\partial \mathbf{f}}{\partial \mathbf{x}}(\mathbf{0}) - \tfrac{\partial \mathcal{N}}{\partial \mathbf{x}}(\mathbf{0})\right] \mathbf{x} + \mathcal{N}(\mathbf{x}) - \mathcal{N}(\mathbf{0})\right) \right\|_1$$
$$\leq \left\| \mathbf{f}(\mathbf{x}) - \left[\tfrac{\partial \mathbf{f}}{\partial \mathbf{x}}(\mathbf{0})\right] \mathbf{x} - \mathbf{g}(\mathbf{x}) \right\|_1$$
$$+ \|\mathbf{g}(\mathbf{x}) - \mathcal{N}(\mathbf{x}) - \mathcal{N}(\mathbf{0})\|_1 + \left\| \left[\tfrac{\partial \mathcal{N}}{\partial \mathbf{x}}(\mathbf{0})\right] \mathbf{x} \right\|_1$$
$$< \tfrac{\varepsilon^*}{2} + \tfrac{\varepsilon^*}{3} + \tfrac{\varepsilon^*}{6} B \leq \tfrac{\varepsilon}{2} + \tfrac{\varepsilon}{3} + \tfrac{\varepsilon}{6} = \varepsilon,$$

for all $\mathbf{x} \in \mathbb{X}$. ∎

### B. PROOF OF THEOREM 2

From Corollary 1 we can find some $\mathbf{g} \in \mathcal{C}^1\left(\mathbb{X}; \mathbb{R}^d\right)$ satisfying $\mathbf{g}(\mathbf{0}) = \mathbf{0}$, $[\partial \mathbf{g}/\partial \mathbf{x}](\mathbf{0}) = \mathbf{0}$, and

$$\max_{\mathbf{x} \in \mathbb{X}} \left\| \mathbf{f}(\mathbf{x}) - \left[\tfrac{\partial \mathbf{f}}{\partial \mathbf{x}}(\mathbf{0})\right] \mathbf{x} - \mathbf{g}(\mathbf{x}) \right\|_1 < \tfrac{\varepsilon}{2}. \quad (36)$$

Applying [42, Exercise 3.23], for all $\mathbf{x} \in \mathbb{X}$ we can decompose $\mathbf{g}(\mathbf{x}) = [\mathbf{h}(\mathbf{x})] \mathbf{x}$, where $\mathbf{h} \in \mathcal{C}\left(\mathbb{X}; \mathbb{R}^{d \times n}\right)$ is given by $\mathbf{h}(\mathbf{x}) = \int_0^1 [\partial \mathbf{g}/\partial \mathbf{x}](s\mathbf{x}) ds$. Further, since $[\partial \mathbf{g}/\partial \mathbf{x}](\mathbf{0}) = \mathbf{0}$ we have $\mathbf{h}(\mathbf{0}) = \int_0^1 [\partial \mathbf{g}/\partial \mathbf{x}](\mathbf{0}) ds = \mathbf{0}$.

Since $\mathbb{X}$ is bounded, there is some $B > 0$ for which $\max_{\mathbf{x} \in \mathbb{X}} \|\mathbf{x}\|_1 \leq B$. Now given $\varepsilon > 0$, by the universal approximation theorem [28] there exists an NN, $\mathcal{N} \in \mathcal{C}^1\left(\mathbb{X}; \mathbb{R}^{d \times n}\right)$, with $\max_{\mathbf{x} \in \mathbb{X}} \|\mathbf{h}(\mathbf{x}) - \mathcal{N}(\mathbf{x})\|_{1,1} < \varepsilon/(4B)$. In particular,

$$\|\mathcal{N}(\mathbf{0})\|_{1,1} = \|\mathbf{0} - \mathcal{N}(\mathbf{0})\|_{1,1} = \|\mathbf{h}(\mathbf{0}) - \mathcal{N}(\mathbf{0})\|_{1,1} < \tfrac{\varepsilon}{4B}.$$

Therefore, for all $\mathbf{x} \in \mathbb{X}$ we have

$$\|\mathbf{g}(\mathbf{x}) - [\mathcal{N}(\mathbf{x}) - \mathcal{N}(\mathbf{0})] \mathbf{x}\|_1$$
$$= \|[\mathbf{h}(\mathbf{x}) - \mathcal{N}(\mathbf{x}) + \mathcal{N}(\mathbf{0})] \mathbf{x}\|_1$$
$$\leq \left( \|\mathbf{h}(\mathbf{x}) - \mathcal{N}(\mathbf{x})\|_{1,1} + \|\mathcal{N}(\mathbf{0})\|_{1,1} \right) \|\mathbf{x}\|_1$$
$$< \left( \tfrac{\varepsilon}{4B} + \tfrac{\varepsilon}{4B} \right) B$$
$$= \tfrac{\varepsilon}{2}. \quad (37)$$

Applying the triangle inequality to (36) and (37) finishes the proof. ∎